\documentclass[twoside]{article}
\usepackage[english]{babel}
\usepackage[cp1251]{inputenc}
\usepackage{amssymb,amsmath,amsthm}
\newcommand{\no}{\noindent}

\newtheorem{de}{Definition}[section]
\newtheorem{theor}{Theorem}[section]
\newtheorem{pr}{Proposition}[section]
\begin{document}

\begin{center}
{\Large \bf{On isomorphism classes and invariants of low dimensional \\
complex filiform Leibniz algebras (part 2)}}
\end{center}
\begin{center}
I.S. Rakhimov$^1$, \ \ S.K. Said Husain$^2$
\end{center}

Institute for Mathematical Research (INSPEM) $\&$ Department of
Mathematics, FS,

UPM, 43400, Serdang, Selangor Darul
Ehsan, (Malaysia)\\

$^1$ isamiddin@science.upm.edu.my $\&$ risamiddin@mail.ru

$^2$ skartini@science.upm.edu.my

\begin{abstract}
\no In \cite{BR1}, \cite{BR2} an elementary method for obtaining a
classification of complex filiform Leibniz algebras based on
invariants was presented. The present paper is an implementation
of this method in low dimensional cases. We give a complete
classification of a subclass of complex filiform Leibniz algebras
obtained from the naturally graded non-Lie filiform Leibniz
algebras. It is known \cite{AO} that this class can be split into
two subclasses. In this paper we shall consider the second one. We
give a hypothetic formula for the adapted number of isomorphism
classes.
\end{abstract}

 \textbf{2000 MSC:} \textit{17A32, 17B30.}

\textbf{Key-Words:} \textit{filiform Leibniz algebra, invariant,
isomorphism.}

\footnote{I.S.R. Corresponding author}

\thispagestyle{empty}

\section{Introduction}

\qquad \ \ Leibniz algebra was first introduced in the early 90's
of the last century by French mathematician J.-L.Loday as a
non-associative algebra with multiplication satisfying the Leibniz
identity:
$$[x[yz]]=[[xy]z] - [[xz]y].$$

This identity is equivalent to the classical Jacobi identity when
the multiplication is skewsymmetric. Leibniz algebras appear to be
related in a natural way to several topics such as differential
geometry, homological algebra, classical algebraic topology ,
algebraic K-theory, loop spaces, noncommutative geometry, quantum
physics etc., as a generalization of the corresponding
applications of Lie algebras to these topics. Most papers concern
to study of homological problems of Leibniz algebras. In \cite{LP}
J.-L.Loday and T.Pirashvili have described the free Leibniz
algebras, paper \cite{MU} by A.A.Mikhalev and U.U.Umirbaev  is
devoted to solution of the non-commutative analogue of the
Jacobian conjecture in the affirmative for free Leibniz algebras ,
in the spirit of the corresponding result of C.Reutenauer
\cite{R}, V.Shpilrain \cite{S} and U.U.Umirbaev \cite{U}. The
problems concerning Cartan subalgebras and solvability were
studied by Sh.A.Ayupov and B.A.Omirov (see, for instance,
\cite{AAO1}). The notion of simple Leibniz algebras was suggested
by S.Abdulkassymova and A.Dzhumadil'daev \cite{AD}, who obtained
some results concerning special cases of simple Leibniz algebras.

Unfortunately, up to now there is no paper including complete
discussion on comparisons the structural theory of Lie and Leibniz
algebras (one means results like Levi-Malcev decomposition
theorem, Lie-Engel theorem, Malcev reduction theorem, the analogue
of Killing form, Dinkin diagrams, root space decompositions, the
Serre presentation, the theory of highest weight representations,
the Weyl character formula and much more).

Papers \cite{AO}, \cite{BR1}, \cite{BR2} \cite{AOR1}, \cite{AOR2},
\cite{GO} and \cite{RK1} can be considered as a certain progress
towards a classification of nilpotent Leibniz algebras.

The present paper is devoted to the classification problem of a
subclass of non-Lie complex filiform Leibniz algebras in low
dimension cases.

\section{Preliminaries}
\begin{de} \emph{An algebra $L$ over a field $K$
is called a {\it Leibniz algebra} if it satisfies the following
Leibniz identity:
$$[x[yz]]=[[xy]z] - [[xz]y].$$}
\end{de}

Let $LB_{n}(K)$ be a subvariety of $Alg_{n}(K)$ consisting of all
$n$-dimensional Leibniz algebras over $K$. It is invariant under
the isomorphic action of $GL_n(K)$ ("transport of structure"). Two
algebras are isomorphic if and only if they belong to the same
orbit under this action. As a subset of $Alg_{n}(K)$ the set
$LB_{n}(K)$ is specified by the system of equations with respect
to structural constants $\gamma_{ij}^{k}$:
$$\sum\limits_{\emph{l}=1}^{\emph{n}}{(\gamma_{\emph{jk}}^{\emph{l}}
\gamma_{\emph{il}}^{\emph{m}}-\gamma_{\emph{ij}}^{\emph{l}}\gamma_{\emph{lk}}^{\emph{m}}+
\gamma_{\emph{ik}}^{\emph{l}}\gamma_{\emph{lj}}^{\emph{m}})}=0.$$

It is easy to see that if the multiplication in Leibniz algebra
happens to be anticommutative then it is a Lie algebra. So Leibniz
algebras are "noncommutative" generalization of Lie algebras. As
for Lie algebras case they are well known and several
classifications have been given, for instance see \cite{GKh}. But
unless simple Lie algebras the classification problem of all Lie
algebras in general remains a big problem. Yu.I.Malcev \cite{Mal}
reduced the classification of solvable Lie algebras to the
classification of nilpotent Lie algebras. Apparently the first
non-trivial classification of some classes of low-dimensional
nilpotent Lie algebra are due to Umlauf. In his thesis \cite{Um}
he presented the redundant list of nilpotent Lie algebras of
dimension less than seven. He gave also the list of nilpotent Lie
algebras of dimension less than ten admitting a so-called adapted
basis (now, the nilpotent Lie algebras with this property are
called \emph{filiform Lie algebras}). It was shown by M.Vergne
\cite{Vr} the importance of filiform Lie algebras in the study of
variety of nilpotent Lie algebra laws. Paper \cite{GJKh} concerns
the classification problem of low dimensional filiform Lie
algebras.

Further if it is not asserted additionally all algebras assumed to
be over the field of complex numbers $\mathbf{C}$.

Let $L$ be a Leibniz algebra. We put: $ L^1=L,\quad
L^{k+1}=[L^k,L],\enskip k\in N.$

\begin{de} \emph{A Leibniz algebra $L$ is said to be {\it
nilpotent} if there exists an integer $s\in N,$ such that $$L^{1}
\supset L^{2} \supset ... \supset L^{s}=\{0\}.$$}
\end{de}

 The smallest
integer $s$ for that $L^{s}=0$ is called {\it the nilindex} of
$L$.

\begin {de} \emph{An $n$-dimensional Leibniz algebra $L$ is
said to be {\it filiform} if $dim L^i =n-i,$ where $2\le i\le n.$}
\end {de}

Let $Leib_n$ denote the class of all $n-$dimensional non-Lie
filiform Leibniz algebras.

There are two sources to get a classification complex filiform
Leibniz algebras. One of them is the naturally graded non-Lie
filiform Leibniz algebras and another one is the naturally graded
filiform Lie algebras. Here we deal with complex filiform Leibniz
algebras obtained from the naturally graded non-Lie filiform
Leibniz algebras.

\begin{theor} \emph{\cite{AO}. Let $L$ be an element of
$Leib_{n+1}.$ Then there exists a basis $\{e_0,e_1, \ldots, e_n\}$
of $L$ such that the structural constants of $L$ on this basis are
defined by one of the following two forms: }

 $ \mbox{\emph{a)} \emph{(The first class)}} \
\ \mu_1^{\overline{\alpha}, \theta }:= \left\{
\begin{array}{lll}
[e_{0}e_{0}]=e_{2},\\

[e_ie_{0}]=e_{i+1}, \qquad \qquad \qquad \qquad \qquad
\qquad \qquad  \ \ 1\leq i\leq n-1 \\

[e_{0}e_{1}]= \alpha_{3} e_{3}+ \alpha_{4} e_{4}+...+\alpha_{n-1}
e_{n-1}+
\theta e_{n}, \\

[e_{j}e_{1}]=\alpha_{3}e_{j+2}+
\alpha_{4}e_{j+3}+...+\alpha_{n+1-j}e_{n}, \qquad 1\leq j\leq n-2
\end{array}
\right. $

$ \mbox{\emph{b)} \emph{(The second class)}} \ \
\mu_2^{\overline{\beta}, \gamma }:= \left\{
\begin{array}{lll}
[\emph{e}_{0}\emph{e}_{0}]=\emph{e}_{2},\\

[e_{i}e_{0}]=e_{i+1}, \qquad \qquad \qquad \qquad \qquad
\qquad \qquad  \ \ 2\leq i\leq n-1 \\

[e_{0}e_{1}]= \beta_{3} e_{3}+ \beta_{4}
e_{4}+...+\beta_{n} e_{n}, \\

[e_{1}e_{1}]=\gamma e_{n},\\

[e_{j}e_{1}]=\beta_{3}e_{j+2}+
\beta_{4}e_{j+3}+...+\beta_{n+1-j}e_{n}, \qquad 2\leq j\leq n-2
\end{array}
\right.$
\end{theor}

Note that the algebras from the first and the second class never
are isomorphic to each other. The classes we denote as
$FLeib_{n+1}$ and $SLeib_{n+1}$, correspondingly. The basis is
said to be "\emph{adapted}".

In this paper we will deal with the second class of algebras of
the above Theorem, the elements of $SLeib_{n+1}$ will be denoted
as $L(\beta_{3},\beta_{4},...,\beta_{n},\gamma)$, meaning that
they are defined by the parameters
$\beta_{3},\beta_{4},...,\beta_{n},\gamma$.

\section{On adapted changes of basis and Isomorphism criterion for $SLeib_{n+1}$}

\qquad \ \ Here we simplify the isomorphic action of $GL_n$
("transport of structure") on $SLeib_n.$ The details of proofs can
be found in \cite{GO}.

\begin{de} \emph{Let $L\in SLeib_{n+1}.$ A basis $\{e_0, e_1,
... , e_n\}$ of $L$ is said to be adapted if its multiplication
table has the form $\mu_2^{\overline{\beta}, \gamma }.$}
\end{de}
Let $L$ be a Leibniz algebra defined on a vector space $V$ and
$\{e_0, e_1, ... , e_n\}$ be the adapted basis of $L.$

\begin{de} \emph{The basis transformation $f\in GL(V)$ is
said to be adapted for the structure of $L$ if the basis
$\{f(e_0), f(e_1), ... , f(e_n)\}$ is adapted.}
\end{de}
The closed subgroup of $GL(V)$ spanned by the adapted
transformations is denoted by $GL_{ad}(V).$

\begin{de} \emph{The following types of basis transformations
of $SLeib_{n+1}$ are said to be elementary:} \\

{\small $\begin{array}{ll}

$\emph{first type}$ \ - \ \sigma(b, n)= \left\{
\begin{array}{ll}
f(e_0) = e_0   \\
f(e_1) = e_1+ be_n, \\
f(e_{i+1})=[f(e_i), f(e_0)], & \ 2 \leq i\leq {n-1}, \\
f(e_2)=[f(e_0), f(e_0)]
\end{array} \right. \\
\end{array}$}\\

{\small $\begin{array}{ll}

$\emph{second type}$ \ - \  \eta(a, k)= \left\{
\begin{array}{ll}
f(e_0) = e_0 + ae_k  \\
f(e_1) = e_1 \\
f(e_{i+1})=[f(e_i), f(e_0)], & \ 2 \leq i\leq {n-1}, \ 2 \leq k\leq n, \\
f(e_2)=[f(e_0), f(e_0)]
\end{array} \right.
\end{array}$} \\

{\small $\begin{array}{ll}

$\emph{third type}$ \ - \  \delta(a, b, d)= \left\{
\begin{array}{ll}
f(e_0) = ae_0 + be_1  \\
f(e_1) = de_1 - \frac{bd \gamma}{a} e_{n-1}, & \ ad\neq 0 \\
f(e_{i+1})=[f(e_i), f(e_0)], & \ 2 \leq i\leq {n-1},  \\
f(e_2)=[f(e_0), f(e_0)]
\end{array} \right.
\end{array}$}\\
\emph{where} $a, b, d \in \mathbf{C}.$
\end{de}

\begin{pr} \emph{If $f$ is an adapted transformation of
$SLeib_{n+1}$ then
$$f=\sigma(b_n, n)\circ \eta(a_n, n)\circ \eta(a_{n-1}, n-2)\circ ...
\circ \eta(a_2, 2)\circ\delta(a_0, a_1, b_1)$$}
\end{pr}

\begin{pr} \emph{Transformations of the form $\sigma(a,
n),$ $\eta(a, n)$ and $\eta(a,k),$ where $2\leq k \leq n-2,$ $a
\in \textbf{C}$ preserve the structural constants of
$SLeib_{n+1}.$}
\end{pr}

Since superposition of adapted transformations is again an adapted
transformation, the proposition above means that the
transformation $\sigma(b_n, n)\circ \eta(a_n, n)\circ
\eta(a_{n-1}, n-2)\circ ... \circ \eta(a_2, 2)$ does not change
the structural constants of $SLeib_{n+1}.$

Thus the action of $GL_{ad}(V)$ on $SLeib_{n+1}$ can be reduced to
the action of elementary transformation of the type three.

Let $\quad R_a^m(x):=[[...[x,\underbrace{ a],a],...,a]
}_{m-times}, \ \mbox{and} \ R_a^0(x):=x.$

Now due to Proposition 3.2 it is easy to see that for
$SLeib_{n+1}$ the adapted change of basis has the form:

$\begin{array}{ll} \left\{\begin{array}{ll}
e_{0}^{'} = \ Ae_0+Be_1 \\
e_{1}^{'} = \ De_1-\frac{BD\gamma}{A}e_{n-1},\\

e_{2}^{'}=A(A+B)e_2+AB(\alpha_3e_3+...+\alpha_{n-1})+B(A\beta_n+B\gamma)e_n,\\

e_{k}^{'}=A(\sum\limits_{i=0}^{k-2}{C_{k-1}^{k-1-i}A^{k-1-i}B^iR_{e_1}^i(e_{k-i})}
+B^{k-1}R_{e_1}^{k-1}(e_0))
\end{array} \right.
\end{array}\\
\mbox{where}, \quad 3\leq k\leq n \quad \mbox{and} \quad
A,D\in\mathbf{C} \quad \mbox{such that} \quad AD\neq 0.$

Now we state the isomorphism criterion with respect to
$GL_{ad}(V).$

We introduce the following series of functions:

$$\begin{array}{lll}
\psi_{t}(y;z)=\psi_{t}(y;z_{3},z_{4},...,z_{n},z_{n+1})=$$ \\

$$z_{t}- \sum \limits_{k=3}^{t-1}(C_{k-1}^{k-2}yz_{t+2-k}+
C_{k-1}^{k-3}y^{2} \sum \limits_{i_{1}=k+2}^{t} z_{t+3-i_{1}}\cdot
z_{i_{1}+1-k}+C_{k-1}^{k-4}y^{3}\sum\limits_{i_{2}=k+3}^{t}
\sum\limits_{i_{1}=k+3}^{i_{2}}z_{t+3-i_{2}} \cdot
z_{i_{2}+3-i_{1}} \cdot z_{i_{1}-k}+...+$$ \\ \newline
$$C_{k-1}^{1}y^{k-2}\sum\limits_{i_{k-3}=2k-2}^{t}
\sum\limits_{i_{k-4}=2k-2}^{i_{k-3}}...
\sum\limits_{i_{1}=2k-2}^{i_{2}}z_{t+3-i_{k-3}} \cdot
z_{i_{k-3}+3-i_{k-4}} \cdot ...\cdot z_{i_{2}+3-i_{1}} \cdot
z_{i_{1}+5-2k}$$ \\ \newline $$+y^{k-1}
\sum\limits_{i_{k-2}=2k-1}^{t}
\sum\limits_{i_{k-3}=2k-1}^{i_{k-2}}...\sum
\limits_{i_{1}=2k-1}^{i_{2}}z_{t+3-i_{k-2}} \cdot
z_{i_{k-2}+3-i_{k-3}} \cdot...\cdot z_{i_{2}+3-i_{1}} \cdot
z_{i_{1}+4-2k})\cdot \psi_{k}(y;z), \end{array}$$ where $3\leq t
\leq n,$

$$\psi_{n+1}(y;z)=z_{n+1}$$

\begin{theor}
\emph{Two algebras $L(\beta)$ and $L(\beta')$
from $SLeib_{n+1}$, where $\beta=(\beta_{3},\beta_{4},...,
\beta_{\emph{n}},\gamma)$, and
$\beta'=(\beta'_{3},\beta'_{4},...,\beta'_{n},\gamma')$, are
isomorphic if and only if there exist complex numbers
$\emph{A},\emph{B}$ and $\emph{D}$ such that
$\emph{A}\emph{D}\neq$ 0 and the following conditions hold:}

\begin{equation}
\beta^{\prime}_{t}=\frac{1}{A^{t-2}}\frac{D}{A}\psi_{t}(\frac{B}{A};\beta),
\end{equation}

$3 \leq t \leq n-1,$

\begin{equation}
\begin{array}{ll}
\beta^{\prime}_{n}=\frac{1}{A^{n-2}}\frac{D}{A}\frac{B}{A}\gamma+\psi_{n}(
\frac{B}{A};\beta),
\end{array}
\end{equation}

\emph{and}

$$\gamma^{\prime}=\frac{1}{A^{n-2}}(\frac{D}{A})^2\psi_{n+1}(\frac{A}{B
};\beta).$$
\end{theor}

To simplify notation let us agree that in the above case for
transition from $(n+1)$-dimensional filiform Leibniz algebra
$L(\beta)$ to $(n+1)$-dimensional filiform Leibniz algebra
$L(\beta^{\prime})$\ \ we write $
\beta^{\prime}=\varrho(\frac{1}{A},\frac{B}{A},\frac{D}{A};\beta)):$

\begin{center}
$\varrho(\frac{1}{A},\frac{B}{A},\frac{D}{A};\beta)=
(\varrho_{1}(\frac{1}{A},\frac{B}{A},\frac{D}{A};\beta),
\varrho_{2}(\frac{1}{A},\frac{B}{A},\frac{D}{A};\beta),
...,\varrho_{n-1}(\frac{1}{A},\frac{B}{A},\frac{D}{A};\beta)),$ \\
\end{center}
where
\begin{center}
$\varrho_{t}(x,y,u;z)=x^{t-1}u\psi_{t+2}(y;z)$ \ \ for \ \ $1 \leq
t \leq n-2,$
\end{center}
and
\begin{center}
 $\varrho_{n-1}(x,y,u;z)=
x^{n-5}u^{2}\psi_{n+1}(y;z).$
\end{center}

Here are the main properties, used in this paper, of the operator
$\varrho$:\newline

$1^{0}. \ \ \varrho(1,0,1;\cdot)) \ \ \mbox{is the identity
operator}.$\newline

$2^{0}. \ \
\varrho(\frac{1}{A_{2}},\frac{B_{2}}{A_{2}},\frac{D_{2}}{A_{2}};
\varrho(\frac{1}{A_{1}},\frac{B_{1}}{A_{1}},
\frac{D_{1}}{A_{1}};\beta))=
\varrho(\frac{1}{A_{1}A_{2}},\frac{B_{1}A_{2}+B_{2}D_{1}}{A_{1}A_{2}},
\frac{D_{1}D_{2}}{A_{1}A_{2}};\beta) $\newline

$3^{0} \ \ \mbox{If} \ \
\beta^{\prime}=\varrho(\frac{1}{A},\frac{B}{A}, \frac{D}{A};\beta)
\ \ \mbox{then} \ \ \beta=\varrho(A,-\frac{B}{D},
\frac{A}{D};\beta^{\prime}).$\\

From here on $n$ is a positive integer. We assume that $n\geq4$
since there are complete classifications of complex nilpotent
Leibniz algebras of dimension at most four \cite{AOR}.

In our investigation we proceed from the viewpoint of \cite{BR2}.
Let $N_{n}$ stands by the adapted number of isomorphism classes in
$SLeib_n.$ (i.e. the number of isomorphism classes on the adapted
basis providing that each parametric family is considered as a one
class).

\section{Classification}

For the simplification purpose we establish the following
notations: let $\Lambda _{1}=4\beta _{3}\beta _{5}-5\beta
_{4}^{2},\ \ \Lambda _{2}=2\beta _{3}^{2}\beta _{6}-6\beta
_{3}\beta _{4}\beta _{5}+\beta _{4}\gamma +4\beta _{4}^{3}, \ \
\Lambda _{3}=4\beta _{3}^{2}\beta _{6}+2\beta _{4}\gamma -7\beta
_{4}^{3}, \ \ \Lambda _{4}=4\beta _{3}^{2}\beta _{6}-7\beta
_{4}^{3}, \ \ \Lambda _{5}=\beta _{3}^{2}\beta _{6}-3\beta
_{3}\beta _{4}\beta _{5}+2\beta _{4}^{3}, \ \ \Lambda _{6}=4\beta
_{3}\beta _{4}\gamma +8\beta _{3}^{3}\beta _{7}-28\beta
_{3}^{2}\beta _{4}\beta _{6}+28\beta _{4}^{4}, \ \ \Lambda
_{7}=4\beta _{3}\beta _{4}\gamma +8\beta _{3}^{3}\beta
_{7}-21\beta _{4}^{4}$ and the letters $\Lambda $ with
\textquotedblleft \ \ ' \textquotedblright denote the same
expression depending on parameters $\beta _{3}^{\prime },\beta
_{4}^{\prime },\beta _{5}^{\prime },\beta _{6}^{\prime }$,$ \beta
_{7}^{\prime }$ and $\gamma ^{\prime }.$

\subsection{Dimension 5}
\qquad \ \  In this section one considers $SLeib_{5}.$ Consider
the following representation of $SLeib_5$ as a disjoint union of
its subsets:

$\qquad SLeib_{5}=U_{1}\bigcup U_{2}\bigcup U_{3}\bigcup
U_{4}\bigcup F,$

where

$\qquad U_{1}=\{L(\beta )\in SLeib_{5}:\beta _{3}\neq 0,\gamma
-2\beta _{3}^{2}\neq 0\},$

$\qquad U_{2}=\{L(\beta )\in SLeib_{5}:\beta _{3}\neq 0,\gamma
-2\beta _{3}^{2}=0,\beta _{4}\neq 0\},$

$\qquad U_{3}=\{L(\beta )\in SLeib_{5}:\beta _{3}=0,\gamma \neq
0\},$

$\qquad U_{4}=\{L(\beta )\in SLeib_{5}:\beta _{3}=0,\gamma
=0,\beta _{4}=0\}, $

$\qquad F=\{L(\beta )\in SLeib_{5}:\beta _{3}=0,\gamma
=0,\beta _{4}\neq 0\}.$\\

\begin{pr}
\emph{Two algebras $L(\beta )$ and $L(\beta ^{\prime })$ from
$U_{1}$ are isomorphic if and only if}
\[
\frac{\gamma }{\beta _{3}^{2}}=\frac{\gamma ^{\prime }}{\beta _{3}^{\prime 2}%
}.
\]
\end{pr}

Thus the algebras from the set $U_{1}$ can be parameterized as $
L(1,0,\lambda).$\\

\begin{pr}\emph{  }

\emph{a) Algebras from  $U_{2}$ are isomorphic to } $L(1,1,2);$

\emph{b) Algebras from $U_{3}$ are isomorphic to} $L(0,0,1);$

\emph{c) Algebras from $U_{4}$ are isomorphic to} $L(0,1,0);$

\emph{d) Algebras from $F$ are isomorphic to} $L(0,0,0).$\\
\end{pr}

\begin{theor}
\emph{Let $L$ be a non-Lie complex filiform Leibniz algebra in
$SLeib_{5}$. Then it is isomorphic to one of the following
pairwise non-isomorphic Leibniz algebras: }

$1)\ \ L(1,0,\lambda ):$

$\qquad \qquad e_{0}e_{0}=e_{2},\ \ e_{i}e_{0}=e_{i+1,}\ \ 2\leq
i\leq 3,\ \ e_{0}e_{1}=e_{3},\ \ e_{1}e_{1}=\lambda  e_{4},\ \
e_{2}e_{1}=e_{4},\ \ $
$\lambda  \in C.$\\

$2)\ \ L(1,1,2):$

$ \qquad \qquad e_{0}e_{0}=e_{2},\ \ e_{i}e_{0}=e_{i+1,}\ \ 2\leq
i\leq 3,\ \ e_{0}e_{1}=e_{3}+e_{4},\ \ e_{1}e_{1}=2e_{4},\ \ e_{2}e_{1}=e_{4}.$\\

$3)\ \ L(0,0,1):$

$\qquad \qquad e_{0}e_{0}=e_{2},\ \ e_{i}e_{0}=e_{i+1,}\ \ 2\leq
i\leq 3,\ \ e_{1}e_{1}=e_{4}. $\\

$ 4)\ \ L(0,1,0):$

$\qquad \qquad e_{0}e_{0}=e_{2},\ \ e_{i}e_{0}=e_{i+1,}\ \ 2\leq
i\leq 3,\ \ e_{0}e_{1}=e_{4}.$\\

$5)\ \ L(0,0,0):$

$\qquad \qquad e_{0}e_{0}=e_{2},\ \ e_{i}e_{0}=e_{i+1,}\ \ 2\leq
i\leq 3.$
\end{theor}

The adapted number of isomorphism classes $N_{5}$=5.

\section{\textbf{Dimension 6}}

 \qquad \ \ This section is devoted to six dimensional case.

 The class $SLeib_6$ can be represented as a disjoint union of several
open and closed subsets:

$\qquad SLeib_{6}=U_{1}\bigcup U_{2}\bigcup U_{3}\bigcup
U_{4}\bigcup U_{5}\bigcup U_{6}\bigcup U_{7}\bigcup U_{8}\bigcup
F,$

where

$\qquad U_{1}=\{L(\beta )\in SLeib_{6}:\beta _{3}\neq 0,\gamma
\neq 0\},$

$\qquad U_{2}=\{L(\beta )\in SLeib_{6}:\beta _{3}\neq 0,\gamma
=0,\Lambda _{1}\neq 0\},$

$\qquad U_{3}=\{L(\beta )\in SLeib_{6}:\beta _{3}\neq 0,\gamma
=0,\Lambda _{1}=0\},$

$\qquad U_{4}=\{L(\beta )\in SLeib_{6}:\beta _{3}=0,\beta _{4}\neq
0,\gamma \neq 0\},$

$\qquad U_{5}=\{L(\beta )\in SLeib_{6}:\beta _{3}=0,\beta _{4}\neq
0,\gamma =0,\beta _{5}\neq 0\},$

$\qquad U_{6}=\{L(\beta )\in SLeib_{6}:\beta _{3}=0,\beta _{4}\neq
0,\gamma =0,\beta _{5}=0\},$

$\qquad U_{7}=\{L(\beta )\in SLeib_{6}:\beta _{3}=0,\beta
_{4}=0,\gamma \neq 0\},$

$\qquad U_{8}=\{L(\alpha )\in SLeib_{6}:\beta _{3}=0,\beta
_{4}=0,\gamma =0,\beta _{5}\neq 0\},$

$\qquad F=\{L(\alpha )\in SLeib_{6}:\beta _{3}=0,\beta
_{4}=0,\gamma =0,\beta _{5}=0\}$.\\

\begin{pr}
 \emph{Two algebras $L(\beta )$ and $L(\beta ^{\prime
})$ from $U_{1}$ are isomorphic if and only if}
\[
\frac{2\beta _{3}\beta _{4}\gamma +\beta _{3}^{2}\Lambda
_{1}}{\gamma ^{2}}=\frac{2\beta _{3}^{\prime }\beta _{4}^{\prime
}\gamma ^{\prime }+\beta _{3}^{\prime 2}\Lambda' _{1}}{\gamma
^{\prime 2}}.
\]
\end{pr}

Thus the set $U_{1}$ can be parameterized as
$L(1,0,\lambda,1).$\\

\begin{pr}\emph{  }

\emph{a) Algebras from $U_{2}$ are isomorphic to} $L(1,0,1,0);$

\emph{b) Algebras from $U_{3}$ are isomorphic to} $L(1,0,0,0);$

\emph{c) Algebras from $U_{4}$ are isomorphic to} $L(0,1,0,1);$

\emph{d) Algebras from $U_{5}$ are isomorphic to} $L(0,1,1,0);$

\emph{e) Algebras from $U_{6}$ are isomorphic to} $L(0,1,0,0);$

\emph{f) Algebras from $U_{7}$ are isomorphic to} $L(0,0,0,1);$

\emph{g) Algebras from $U_{8}$ are isomorphic to} $L(0,0,1,0);$

\emph{h) Algebras from $F$ are isomorphic to} $L(0,0,0,0).$
\end{pr}

\begin{theor}

\emph{Let $L$ be a non-Lie complex filiform Leibniz algebra in
$SLeib_{6}$. Then it is isomorphic to one of the following
pairwise non-isomorphic Leibniz algebras: }

 $1)\ \ L(1,0,\lambda
,1):$

$ \qquad \qquad e_{0}e_{0}=e_{2},\ \ e_{i}e_{0}=e_{i+1,}\ \ 2\leq
i\leq 4,\ \ e_{0}e_{1}=e_{3}+\lambda  e_{5},\ \ e_{1}e_{1}=e_{5},\
\ e_{2}e_{1}=e_{4},$

$ \qquad \qquad
e_{3}e_{1}=e_{5},\ \ \lambda  \in C.$\\

$ 2)\ \ L(1,0,1,0):$

$\qquad \qquad e_{0}e_{0}=e_{2},\ \ e_{i}e_{0}=e_{i+1,}\ \ 2\leq
i\leq 4,\ \ e_{0}e_{1}=e_{3}+e_{5},\ \ e_{2}e_{1}=e_{4},\ \ e_{3}e_{1}=e_{5}.$\\

$ 3)\ \ L(1,0,0,0):$

$\qquad \qquad e_{0}e_{0}=e_{2},\ \ e_{i}e_{0}=e_{i+1,}\ \ 2\leq
i\leq 4,\ \ e_{0}e_{1}=e_{3}, \ \ e_{2}e_{1}=e_{4},\ \ e_{3}e_{1}=e_{5}.$\\

$ 4)\ \ L(0,1,0,1):$

$\qquad \qquad e_{0}e_{0}=e_{2},\ \ e_{i}e_{0}=e_{i+1,}\ \ 2\leq
i\leq 4,\ \
e_{0}e_{1}=e_{4}, \ \ e_{1}e_{1}=e_{5},\ \ e_{2}e_{1}=e_{5}.$\\

$ 5)\ \ L(0,1,1,0)$

$\qquad \qquad e_{0}e_{0}=e_{2},\ \ e_{i}e_{0}=e_{i+1,}\ \ 2\leq
i\leq 4,\ \ e_{0}e_{1}=e_{4}+e_{5},\ \ e_{2}e_{1}=e_{5}.$\\

$ 6) \ \ L(0,1,0,0):$

$ \qquad \qquad e_{0}e_{0}=e_{2},\ \ e_{i}e_{0}=e_{i+1,}\ \ 2\leq
i\leq 4,\ \ e_{0}e_{1}=e_{4}, \ \ e_{2}e_{1}=e_{5}.$\\

$ 7)\ L(0,0,0,1):$

 $\qquad \qquad  e_{0}e_{0}=e_{2},\ \ e_{i}e_{0}=e_{i+1,}\ \ 2\leq
i\leq 4,\ \ e_{1}e_{1}=e_{5}.$\\

$ 8)\ \ L(0,0,1,0):$

$ \qquad \qquad e_{0}e_{0}=e_{2},\ \ e_{i}e_{0}=e_{i+1,}\ \ 2\leq
i\leq 4,\ \ e_{0}e_{1}=e_{5}.$\\

$9)\ \ L(0,0,0,0):$

$\qquad \qquad e_{0}e_{0}=e_{2},\ \ e_{i}e_{0}=e_{i+1,}\ \ 2\leq
i\leq 4.$
\end{theor}

The adapted number of isomorphism classes $N_{6}$=9.

\section{\textbf{Dimension 7}}

Consider the following decomposition of $SLeib_{7}$ into its
disjoint subsets:

\qquad \ \ $SLeib_{7}=U_{1}\bigcup U_{2}\bigcup U_{3}\bigcup
U_{4}\bigcup U_{5}\bigcup U_{6}\bigcup U_{7}\bigcup U_{8}\bigcup
U_{9}\bigcup $ $U_{10}\bigcup U_{11}\bigcup U_{12}\bigcup
U_{13}\bigcup U_{14}\bigcup F,$

where

$U_{1}=\{L(\beta )\in SLeib_{7}:\beta _{3}\neq 0,\Lambda _{1}\neq
0, \Lambda _{2}\neq 0\},$

$ U_{2}=\{L(\beta )\in SLeib_{7}:\beta _{3}\neq 0,\Lambda _{1}\neq
0, \Lambda _{2}=0,\gamma \neq 0\},$

$ U_{3}=\{L(\beta )\in SLeib_{7}:\beta _{3}\neq 0,\Lambda
_{1}=0,\Lambda _{3}\neq 0,\gamma \neq 0\},$

$U_{4}=\{L(\beta )\in SLeib_{7}:\beta _{3}=0,\beta _{4}\neq
0,\beta _{5}\neq 0 \},$

$U_{5}=\{L(\beta )\in SLeib_{7}:\beta _{3}=0,\beta _{4}\neq
0,\beta _{5}= 0,\gamma-3\beta _{4}^{2}\neq 0 \},$

$ U_{6}=\{L(\beta )\in SLeib_{7}:\beta _{3}\neq 0,\Lambda
_{1}=0,\Lambda _{3}= 0,\gamma = 0,\Lambda _{4}\neq 0\},$

$U_{7}=\{L(\beta )\in SLeib_{7}:\beta _{3}=0,\beta _{4}\neq
0,\beta _{5}= 0,\gamma-3\beta _{4}^{2}= 0 \},$

$U_{8}=\{L(\beta )\in SLeib_{7}:\beta _{3}=0,\beta _{4}=0,\beta
_{5}\neq 0,\beta _{6}\neq 0,\gamma \neq 0 \},$

$U_{9}=\{L(\beta )\in SLeib_{7}:\beta _{3}=0,\beta _{4}=0,\beta
_{5}\neq 0,\beta _{6}\neq 0,\gamma = 0\},$

 $U_{10}=\{L(\beta )\in SLeib_{7}:\beta _{3}=0,\beta _{4}=0,\beta
_{5}\neq 0,\beta _{6}=0,\gamma \neq 0\},$

 $U_{11}=\{L(\beta )\in SLeib_{7}:\beta _{3}=0,\beta _{4}=0,\beta
_{5}\neq 0,\beta _{6}=0,\gamma =0\},$

$U_{12}=\{L(\beta )\in SLeib_{7}:\beta _{3}=0,\beta _{4}=0,\beta
_{5}=0,\beta _{6}\neq 0,\gamma \neq 0\},$

$U_{13}=\{L(\beta )\in SLeib_{7}:\beta _{3}=0,\beta _{4}=0,\beta
_{5}=0,\beta _{6}\neq 0,\gamma =0\},$

$U_{14}=\{L(\beta )\in SLeib_{7}:\beta _{3}=0,\beta _{4}=0,\beta
_{5}=0,\beta _{6}=0,\gamma \neq 0\},$

$F=\{L(\beta )\in SLeib_{7}:\beta _{3}=0,\beta _{4}=0,\beta
_{5}=0,\beta _{6}=0,\gamma =0\}.$\\

 \begin{pr}
 \emph{Two algebras $L(\beta )$ and $L(\beta ^{\prime })$ from
$U_{1}$ are isomorphic if and only if}
\[
\frac{\Lambda _{1}^{3}}{\Lambda _{2}^{2}\ }=\frac{\Lambda _{1}^{\prime 3}}{%
\Lambda _{2}^{\prime 2}\ }
\]
\[
\frac{\gamma \Lambda _{1}^{2}}{\Lambda _{2}^{2}\ }=\frac{\gamma
^{\prime }\Lambda _{1}^{\prime 2}}{\Lambda _{2}^{\prime 2}\ }
\]
\end{pr}

Thus the algebras from the set $U_{1}$ can be parameterized as
$L(1,0,\lambda _{1},\lambda _{1},\lambda _{2}).$\\

 \begin{pr}
 \emph{ Two algebras $L(\beta )$ and
$L(\beta ^{\prime })$ from $U_{2}$ are isomorphic if and only if}
\[
\frac{\gamma}{\Lambda_1 }=\frac{\gamma'}{\Lambda'_1 }.
\]
\end{pr}

Thus the algebras from the set $U_{2}$ can be parameterized as
$L(1,0,1,0,\lambda ).$\\

 \begin{pr}
 \emph{Two algebras $L(\beta )$ and
$L(\beta ^{\prime })$ from $U_{3}$ are isomorphic if and only if}
\[
\frac{\gamma ^{3}}{\Lambda_1 ^{2}}=\frac{\gamma ^{\prime
3}}{\Lambda_{1}^{\prime 2}}.
\]
\end{pr}

Thus the algebras from the set $U_{3}$ can be parameterized as
$L(1,0,0,\lambda,\lambda ).$\\

 \begin{pr}
 \emph{ Two algebras $L(\beta )$ and
$L(\beta ^{\prime })$ from $U_{4}$ are isomorphic if and only if}
\[
\frac{\gamma }{\beta _{4}^{2}}=\frac{\gamma ^{\prime }}{\beta _{4}^{\prime 2}%
}.
\]
\end{pr}

Thus the algebras from the set $U_{4}$ can be parameterized as
$L(0,1,1,0,\lambda ).$\\

\begin{pr}
\emph{Two algebras $L(\beta )$ and $L(\beta ^{\prime })$ from
$U_{5}$ are isomorphic if and only if}
\[
\frac{\gamma }{\beta _{4}^{2}}=\frac{\gamma ^{\prime }}{\beta _{4}^{\prime 2}%
}.
\]
\end{pr}

Thus the algebras from the set $U_{5}$ can be parameterized as
$L(0,1,0,0,\lambda ).$\\

\begin{pr}\emph{  }

\emph{a) Algebras from $U_{6}$ are isomorphic to} $L(1,0,0,1,0);$

\emph{b) Algebras from $U_{7}$ are isomorphic to} $L(0,1,0,1,3);$

\emph{c) Algebras from $U_{8}$ are isomorphic to} $L(0,0,1,0,1);$

\emph{d) Algebras from $U_{9}$ are isomorphic to} $L(0,0,1,1,0);$

\emph{e) Algebras from $U_{10}$ are isomorphic to} $L(0,0,1,1,1);$

\emph{f) Algebras from $U_{11}$ are isomorphic to} $L(0,0,1,0,0);$

\emph{g) Algebras from $U_{12}$ are isomorphic to} $L(0,0,0,0,1);$

\emph{h) Algebras from $U_{13}$ are isomorphic to} $L(0,0,0,1,0);$

\emph{i) Algebras from $U_{14}$ are isomorphic to} $L(0,0,0,1,1);$

\emph{j) Algebras from $F$ are isomorphic to} $L(0,0,0,0,0).$
\end{pr}

 \begin{theor}

\emph{ Let $L$ be non-Lie complex filiform Leibniz algebra in
$SLeib_{7}$. Then it is isomorphic to one of the following
pairwise non-isomorphic Leibniz algebras:}

$1)\ \ L(1,0,\lambda_{1},\lambda_{1},\lambda_{2}):$

$ \qquad \qquad e_{0}e_{0}=e_{2},\ \ e_{i}e_{0}=e_{i+1,}\ \ 2\leq
i\leq 5,\ \ e_{0}e_{1}=e_{3}+\lambda_{1}e_{5}+\lambda _{1}e_{6},\
\ e_{1}e_{1}=\lambda_{2}e_{6},$

$ \qquad \qquad e_{2}e_{1}=e_{4}+\lambda_{1}e_{6},\ \
e_{3}e_{1}=e_{5},\ \ e_{4}e_{1}=e_{6},\ \ \lambda
_{1},\lambda_{2}\in C.$\\

 $2)\ \  L(1,0,1, 0,\lambda):$

$\qquad \qquad e_{0}e_{0}=e_{2},\ \ e_{i}e_{0}=e_{i+1,}\ \ 2\leq
i\leq 5,\ \ e_{0}e_{1}=e_{3}+ e_{5},\ \ e_{1}e_{1}=\lambda e_{6},\
\ e_{2}e_{1}=e_{4}+ e_{6},$

$ \qquad \qquad  e_{3}e_{1}=e_{5},\ \ e_{4}e_{1}=e_{6},\ \ \lambda \in C.$\\

$3)\ \  L(1,0,0,\lambda ,\lambda ):$

$\qquad \qquad e_{0}e_{0}=e_{2},\ \ e_{i}e_{0}=e_{i+1,}\ \ 2\leq
i\leq 5,\ \ e_{0}e_{1}=e_{3}+\lambda e_{6},\ \ e_{1}e_{1}=\lambda
e_{6},\ \ e_{2}e_{1}=e_{4},$

$ \qquad \qquad  e_{3}e_{1}=e_{5},\ \ e_{4}e_{1}=e_{6},\ \ $
$\lambda \in C.$\\

$4)\ \ L(0,1,1,0,\lambda  ):$

 $ \qquad \qquad e_{0}e_{0}=e_{2},\ \ e_{i}e_{0}=e_{i+1,}\ \ 2\leq
i\leq 5,\ \ e_{0}e_{1}=e_{4}+e_{5},\ \ e_{1}e_{1}=\lambda e_{6},\
\ e_{2}e_{1}=e_{5}+e_{6},$

$ \qquad \qquad  e_{3}e_{1}=e_{6},\ \ \lambda \in
C.$\\

$5)\ \ L(0,1,0,0,\lambda  ):$

$\qquad \qquad e_{0}e_{0}=e_{2},\ \ e_{i}e_{0}=e_{i+1,}\ \ 2\leq
i\leq 5,\ \ e_{0}e_{1}=e_{4}+e_{5}+e_{6},\ \ e_{1}e_{1}=\lambda
e_{6},$

$ \qquad \qquad  e_{2}e_{1}=e_{5}+e_{6},\ \ e_{3}e_{1}=e_{6},\ \
\lambda \in
C.$\\

$6)\ \ L(1,0,0,1,0 ):$

 $ \qquad \qquad e_{0}e_{0}=e_{2},\ \ e_{i}e_{0}=e_{i+1,}\ \ 2\leq
i\leq 5,\ \ e_{0}e_{1}=e_{3}+e_{4},\ \ e_{2}e_{1}=e_{4},\ \
e_{3}e_{1}=e_{6},$

$ \qquad \qquad e_{4}e_{1}=e_{6}.$\\

$7)\ \ L(0,1,0,1,3 ):$

$ \qquad \qquad e_{0}e_{0}=e_{2},\ \ e_{i}e_{0}=e_{i+1,}\ \ 2\leq
i\leq 5,\ \  e_{0}e_{1}=e_{4}+e_{6},\ \ e_{1}e_{1}=3 e_{6},\ \
e_{2}e_{1}=e_{5},$

$ \qquad \qquad  e_{3}e_{1}=e_{6}.$\\

$8)\ \  L(0,0,1,0,1 ):$

$ \qquad \qquad e_{0}e_{0}=e_{2},\ \ e_{i}e_{0}=e_{i+1,}\ \ 2\leq
i\leq 5,\ \ e_{0}e_{1}=e_{5},\ \ e_{1}e_{1}=e_{6},\ \
e_{2}e_{1}=e_{6}.$\\

$9)\ \ L(0,0,1,1,0):$

 $\qquad \qquad e_{0}e_{0}=e_{2},\ \ e_{i}e_{0}=e_{i+1,}\ \ 2\leq
i\leq 5,\ \ e_{0}e_{1}=e_{5}+e_{6},\ \ e_{2}e_{1}=e_{6}.$\\

$10)\ \  L(0,0,1,1,1):$

 $ \qquad \qquad e_{0}e_{0}=e_{2},\ \ e_{i}e_{0}=e_{i+1,}\ \ 2\leq
i\leq 5,\ \ e_{0}e_{1}=e_{5}+e_{6},\ \ e_{1}e_{1}=e_{6},\ \
e_{2}e_{1}=e_{6}.$\\

$11)\ \  L(0,0,1,0,0):$

 $ \qquad \qquad e_{0}e_{0}=e_{2},\ \ e_{i}e_{0}=e_{i+1,}\ \ 2\leq
i\leq 5,\ \ e_{0}e_{1}=e_{5},\ \ e_{2}e_{1}=e_{6}.$\\

 $12)\ \ L(0,0,0,0,1):$

 $ \qquad \qquad e_{0}e_{0}=e_{2},\ \ e_{i}e_{0}=e_{i+1,}\ \ 2\leq
i\leq 5,\ \  e_{1}e_{1}=e_{6}.$\\

$13)\ \ L(0,0,0,1,0):$

 $ \qquad \qquad e_{0}e_{0}=e_{2},\ \ e_{i}e_{0}=e_{i+1,}\ \ 2\leq
i\leq 5,\ \ e_{0}e_{1}=e_{6}.$\\

$14)\ \  L(0,0,0,1,1):$

$\qquad \qquad e_{0}e_{0}=e_{2},\ \ e_{i}e_{0}=e_{i+1,}\ \ 2\leq
i\leq 5,\ \ e_{0}e_{1}=e_{6},\ \ e_{1}e_{1}=e_{6}.$\\

$15)\ \  L(0,0,0,0,0):$

 $ \qquad \qquad e_{0}e_{0}=e_{2},\ \ e_{i}e_{0}=e_{i+1,}\ \ 2\leq
i\leq 5. $
\end{theor}

\section{\textbf{Dimension 8}}

\qquad \ \ $SLeib_{8}=U_{1}\bigcup U_{2}\bigcup U_{3}\bigcup
U_{4}\bigcup U_{5}\bigcup U_{6}\bigcup U_{7}\bigcup U_{8}\bigcup
U_{9}\bigcup U_{10}\bigcup U_{11}\bigcup  U_{12}\bigcup
U_{13}\bigcup $

$\qquad \qquad \ U_{14}\bigcup U_{15}\bigcup U_{16}\bigcup
U_{17}\bigcup U_{18}\bigcup U_{19}\bigcup U_{20}\bigcup
U_{21}\bigcup U_{22}\bigcup F,$

where

$\ U_{1}=\{L(\beta )\in SLeib_{8}:\beta _{3}\neq 0,\Lambda
_{1}\neq 0, \Lambda _{5}\neq 0\},$

$\ U_{2}=\{L(\beta )\in SLeib_{8}:\beta _{3}\neq 0,\Lambda
_{1}\neq 0, \Lambda _{5}=0,\gamma \neq 0\},$

$\ U_{3}=\{L(\beta )\in SLeib_{8}:\beta _{3}\neq 0,\Lambda
_{1}\neq 0, \Lambda _{5}=0,\gamma =0\},$

$\ U_{4}=\{L(\beta )\in SLeib_{8}:\beta _{3}\neq 0,\Lambda _{1}=
0, \Lambda _{4}\neq 0,\Lambda _{6}\neq 0\},$

$\ U_{5}=\{L(\beta )\in SLeib_{8}:\beta _{3}\neq 0,\Lambda _{1}=
0, \Lambda _{4}\neq 0,\Lambda _{6}=0\},$

\ $U_{6}=\{L(\beta )\in SLeib_{8}:\beta _{3}\neq 0,\Lambda _{1}=
0, \Lambda _{4}=0,\gamma \neq 0,\Lambda _{7}\neq 0\},$

\ $U_{7}=\{L(\beta )\in SLeib_{8}:\beta _{3}=0,\beta _{4}\neq
0,\beta _{5}\neq 0,\beta _{6}\neq 0\},$

\ $U_{8}=\{L(\beta )\in SLeib_{8}:\beta _{3}=0,\beta _{4}\neq
0,\beta _{5}\neq 0,\beta _{6}=0\},$

\ $U_{9}=\{L(\beta )\in SLeib_{8}:\beta _{3}=0,\beta _{4}\neq
0,\beta _{5}=0,\beta _{6}\neq 0,\gamma \neq 0\},$

\ $U_{10}=\{L(\beta )\in SLeib_{8}:\beta _{3}=0,\beta _{4}\neq
0,\beta _{5}=0,\beta _{6}=0,\gamma \neq 0\},$

$\ U_{11}=\{L(\beta )\in SLeib_{8}:\beta _{3}=0,\beta _{4}=0,\beta
_{5}\neq 0,\beta _{6}\neq 0\},$

\ $U_{12}=\{L(\beta )\in SLeib_{8}:\beta _{3}=0,\beta _{4}=0,\beta
_{5}\neq 0,\beta _{6}=0,\gamma =0,\beta _{7}\neq 0\},$

\ $U_{13}=\{L(\beta )\in SLeib_{8}:\beta _{3}=0,\beta _{4}=0,\beta
_{5}=0,\beta _{6}\neq 0,\gamma \neq 0,\beta _{7}\neq 0\},$

\ $U_{14}=\{L(\beta )\in SLeib_{8}:\beta _{3}=0,\beta _{4}=0,\beta
_{5}=0,\beta _{6}\neq 0,\gamma \neq 0,\beta _{7}=0\},$

\ $U_{15}=\{L(\beta )\in SLeib_{8}:\beta _{3}=0,\beta _{4}\neq
0,\beta _{5}=0,\beta _{6}\neq 0,\gamma =0,\beta _{7}\neq 0\},$

\ $U_{16}=\{L(\beta )\in SLeib_{8}:\beta _{3}=0,\beta _{4}\neq
0,\beta _{5}=0,\beta _{6}\neq 0,\gamma =0,\beta _{7}=0\},$

\ $U_{17}=\{L(\beta )\in SLeib_{8}:\beta _{3}=0,\beta _{4}\neq
0,\beta _{5}=0,\beta _{6}=0,\gamma =0,\beta _{7}\neq 0\},$

\ $U_{18}=\{L(\beta )\in SLeib_{8}:\beta _{3}=0,\beta _{4}\neq
0,\beta _{5}=0,\beta _{6}=0,\gamma =0,\beta _{7}=0\},$

\ $U_{19}=\{L(\beta )\in SLeib_{8}:\beta _{3}=0,\beta _{4}=0,\beta
_{5}\neq 0,\beta _{6}=0,\gamma \neq 0\},$

\ $U_{20}=\{L(\beta )\in SLeib_{8}:\beta _{3}=0,\beta _{4}=0,\beta
_{5}=0,\beta _{6}=0,\beta _{7}\neq 0,\gamma \neq 0\},$

$\ U_{21}=\{L(\beta )\in SLeib_{8}:\beta _{3}=0,\beta _{4}=0,\beta
_{5}=0,\beta _{6}=0,\beta _{7}\neq 0,\gamma =0\},$

$\ U_{22}=\{L(\beta )\in SLeib_{8}:\beta _{3}=0,\beta _{4}=0,\beta
_{5}=0,\beta _{6}=0,\beta _{7}=0,\gamma \neq 0\},$

$\ F=\{L(\beta )\in SLeib_{8}:\beta _{3}=0,\beta
_{4}=0,\beta _{5}=0,\beta _{6}=0,\beta _{7}=0,\gamma =0\}.$\\

\begin{pr}

\emph{ Two algebras $L(\beta )$ and $L(\beta ^{\prime })$ from
$U_{1}$ are isomorphic if and only if}
\[
\frac{\Lambda_1 ^{3}}{\Lambda_5 ^{2}}=\frac{ \Lambda_1 ^{\prime
3}}{ \Lambda_5 ^{\prime 2}}
\]

\begin{eqnarray*}
&&\frac{\Lambda_1 ^{4}\left(\Lambda_7 -28\beta_{4} \Lambda_5
-14\beta_{4}^{2} \Lambda_1\right) }{\Lambda_5 ^{4}}
=\frac{\Lambda_1 ^{\prime 4}\left(  \Lambda_7 ^{\prime}
-28\beta_{4} ^{\prime} \Lambda_5 ^{\prime} -14\beta_{4}^{\prime 2}
\Lambda_1 ^{\prime} \right)}{\Lambda_5 ^{\prime 4}}
\end{eqnarray*}

\[
\frac{\beta _{3}\gamma \Lambda_1^{3}}{ \Lambda_5 ^{3}
}=\frac{\beta _{3}^{\prime}\gamma ^{\prime}\Lambda_1 ^{\prime 3}}{%
 \Lambda_5 ^{\prime 3}}.
\]
\end{pr}
Thus the algebras from the set $U_{1}$ can be parameterized as
$L(1,0,\lambda_{1},\lambda _{1},\lambda _{2},\lambda _{3}).$\\

\begin{pr}
\emph{Two algebras $L(\beta )$ and $L(\beta ^{\prime })$ from
$U_{2}$ are isomorphic if and only if}
\[
\frac{ \Lambda_1 ^{3}}{\beta _{3}^{2}\gamma ^{2}}=\frac{ \Lambda_1
^{\prime 3}}{\beta _{3}^{\prime 2}\gamma ^{\prime 2}}
\]%
\begin{eqnarray*}
&&\frac{ \Lambda_1
^{4}\left(\Lambda_7-14\beta_{4}^{2}\Lambda_1\right) }{\beta
_{3}^{4}\gamma ^{4}} = \frac{ \Lambda_1 ^{\prime 4}\left(
\Lambda_7 ^{\prime} -14\beta_{4}^{\prime 2} \Lambda_1 ^{\prime}
\right)}{\beta _{3}^{\prime 4}\gamma ^{\prime 4}}.
\end{eqnarray*}

\end{pr}
Thus the algebras from the set $U_{2}$ can be parameterized as
$L(1,0,\lambda_{1},0,\lambda _{2},\lambda _{1}).$\\

 \begin{pr}
\emph{Two algebras $L(\beta )$ and $L(\beta ^{\prime })$ from
$U_{3}$ are isomorphic if and only if}
\[
\frac{\Lambda_1 ^{4}\left(
\Lambda_7-14\beta_{4}^{2}\Lambda_1-4\beta_{3}\beta_{4}\gamma\right)
}{\Lambda_1^{2}}=\frac{\Lambda_1 ^{\prime 4}\left( \Lambda_7
^{\prime} -14\beta_{4}^{\prime 2}
\Lambda_1^{\prime}-4\beta_{3}^{\prime}\beta_{4}^{\prime}\gamma^{\prime}\right)}{\Lambda_1
^{\prime 2}}.
\]
\end{pr}
Thus the algebras from the set $U_{3}$ can be parameterized as
$L(1,0,1,0,\lambda ,0).$\\

\begin{pr}
\emph{Two algebras $L(\beta )$ and $L(\beta ^{\prime })$ from
$U_{4}$ are isomorphic if and only if}
\[
\frac{ \Lambda_4 ^{4}}{ \Lambda_6 ^{3}}=\frac{ \Lambda_4 ^{\prime 4}}{%
 \Lambda_6 ^{\prime 3}}
\]
\[
\frac{\beta _{3}\gamma  \Lambda_4^{3}}{ \Lambda_6
^{3}}=\frac{\beta _{3}^{\prime }\gamma ^{\prime } \Lambda_4 ^{\prime 3}%
}{ \Lambda_6 ^{\prime 3}}.
\]
\end{pr}
Thus the algebras from the set $U_{4}$ can be parameterized as
$L(1,0,0,\lambda _{1},\lambda_{1},\lambda _{2}).$\\

\begin{pr}
\emph{Two algebras $L(\beta )$ and $L(\beta ^{\prime })$ from
$U_{5}$ are isomorphic if and only if}
\[
\frac{\beta _{3}\gamma }{ \Lambda_4 }=\frac{\beta _{3}^{\prime
}\gamma ^{\prime }}{ \Lambda_4 ^{\prime}}.
\]
\end{pr}
Thus the algebras from the set $U_{5}$ can be parameterized as
$L(1,0,0,1,0,\lambda ).$\\

\begin{pr}
\emph{Two algebras $L(\beta )$ and $L(\beta ^{\prime })$ from
$U_{6}$ are isomorphic if and only if}
\[
\frac{\beta _{3}^{4}\gamma ^{4}}{ \Lambda_7 ^{3}}=\frac{\beta
_{3}^{\prime 4}\gamma ^{\prime 4}}{ \Lambda_7 ^{\prime 3}}.
\]
\end{pr}
Thus the algebras from the set $U_{6}$ can be parameterized as
$L(1,0,0,0,\lambda ,\lambda).$\\

\begin{pr}
\emph{Two algebras $L(\beta )$ and $L(\beta ^{\prime })$ from
$U_{7}$ are isomorphic if and only if}
\[
\frac{\beta _{5}\beta _{6}\gamma +3\beta _{4}^{2}\beta _{7}-7\beta
_{4}\beta _{5}^{2}\beta _{6}}{\beta _{5}^{3}}=\frac{\beta
_{5}^{\prime }\beta _{6}^{\prime }\gamma ^{\prime }+3\beta
_{4}^{\prime 2}\beta _{7}^{\prime }-7\beta _{4}^{\prime }\beta
_{5}^{\prime 2}\beta _{6}^{\prime }}{\beta _{5}^{\prime 3}}
\]%
\[
\frac{\gamma }{\beta _{4}\beta _{5}}=\frac{\gamma ^{\prime
}}{\beta _{4}^{\prime }\beta _{5}^{\prime }} .\]
\end{pr}
Thus the algebras from the set $U_{7}$ can be parameterized as
$L(0,1,1,0,\lambda_{1},\lambda_{2}).$\\

\begin{pr}
\emph{Two algebras $L(\beta )$ and $L(\beta ^{\prime })$ from
$U_{8}$ are isomorphic if and only if}
\[
\frac{\beta _{4}\beta _{5}\gamma +3\beta _{4}^{2}\beta _{7}+7\beta
_{4}^{2}\beta _{5}^{2}}{\beta _{5}^{3}}=\frac{\beta _{4}^{\prime
}\beta _{5}^{\prime }\gamma ^{\prime }+3\beta _{4}^{\prime 2}\beta
_{7}^{\prime }+7\beta _{4}^{\prime 2}\beta _{5}^{\prime 2}}{\beta
_{5}^{\prime 3}}
\]
\[
\frac{\gamma }{\beta _{4}\beta _{5}}=\frac{\gamma ^{\prime
}}{\beta _{4}^{\prime }\beta _{5}^{\prime }} .\]
\end{pr}
Thus the algebras from the set $U_{8}$ can be parameterized as
$L(0,1,1,-1,\lambda_{1},\lambda_{2}).$\\

\begin{pr}
\emph{Two algebras $L(\beta )$ and $L(\beta ^{\prime })$ from
$U_{9}$ are isomorphic if and only if}
\[
\frac{\beta _{4}^{3}\left( \beta _{6}\gamma +3\beta _{4}^{2}\beta
_{7}\right) }{\gamma ^{3}}=\frac{\beta _{4}^{\prime 3}\left( \beta
_{6}^{\prime }\gamma ^{\prime }+3\beta _{4}^{\prime 2}\beta
_{7}^{\prime }\right) }{\gamma ^{\prime 3}}
.\]%
\end{pr}
Thus the algebras from the set $U_{9}$ can be parameterized as
$L(0,1,0,0,\lambda,1).$\\

\begin{pr}
\emph{Two algebras $L(\beta )$ and $L(\beta ^{\prime })$ from
$U_{10}$ are isomorphic if and only if}
\[
\frac{\beta _{4}^{5}\beta _{7}+\gamma ^{3}}{\gamma
^{3}}=\frac{\beta
_{4}^{\prime 5}\beta _{7}^{\prime }+\gamma ^{\prime 3}}{\gamma ^{\prime 3}}%
.
\]
\end{pr}
Thus the algebras from the set $U_{10}$ can be parameterized as
$L(0,1,0,-1,\lambda,1).$\\

\begin{pr}
\emph{Two algebras $L(\beta )$ and $L(\beta ^{\prime })$ from
$U_{11}$ are isomorphic if and only if}
\[
\frac{\beta _{6}\gamma }{\beta _{5}^{3}}\ \ \ =\frac{\beta
_{6}^{\prime }\gamma ^{\prime }}{\beta _{5}^{\prime 3}} .
\]
\end{pr}
Thus the algebras from the set $U_{11}$ can be parameterized as
$L(0,0,1,1,0,\lambda).$\\

\begin{pr}
\emph{Two algebras $L(\beta )$ and $L(\beta ^{\prime })$ from
$U_{12}$ are isomorphic if and only if}
\[
\frac{\beta _{7}^{3}}{\beta _{5}^{5}}=\frac{\beta _{7}^{\prime
3}}{\beta _{5}^{\prime 5}} .\]
\end{pr}
Thus the algebras from the set $U_{12}$ can be parameterized as
$L(0,0,1,0,\lambda ,0).$\\

\begin{pr}
\emph{Two algebras $L(\beta )$ and $L(\beta ^{\prime })$ from
$U_{13}$ are isomorphic if and only if}
\[
\frac{\beta _{6}}{\gamma ^{2}}=\frac{\beta _{6}^{\prime }}{\gamma ^{\prime 2}%
}.
\]
\end{pr}
Thus the algebras from the set $U_{13}$ can be parameterized as
$L(0,0,0,1,0,\lambda).$\\

\begin{pr}
\emph{Two algebras $L(\beta )$ and $L(\beta ^{\prime })$ from
$U_{14}$ are isomorphic if and only if}
\[
\frac{\beta _{6}}{\gamma ^{2}}=\frac{\beta _{6}^{\prime }}{\gamma ^{\prime 2}%
}.
\]
\end{pr}
Thus the algebras from the set $U_{14}$ can be parameterized as
$L(0,0,0,1,1,\lambda ).$\\

\begin{pr}\emph{  }

\emph{a) Algebras from $U_{15}$ are isomorphic to}
$L(0,1,0,0,1,0);$

\emph{b) Algebras from $U_{16}$ are isomorphic to}
$L(0,1,0,0,0,0);$

\emph{c) Algebras from $U_{17}$ are isomorphic to}
$L(0,1,0,-1,1,0);$

\emph{d) Algebras from $U_{18}$ are isomorphic to} $%
L(0,1,0,-1,0,0);$

\emph{e) Algebras from $U_{19}$ are isomorphic to}
$L(0,0,1,0,0,1);$

\emph{f) Algebras from $U_{20}$ are isomorphic to}
$L(0,0,0,0,0,1);$

\emph{g) Algebras from $U_{21}$ are isomorphic to}
$L(0,0,0,0,1,0);$

\emph{h) Algebras from $U_{22}$ are isomorphic to}
$L(0,0,0,0,1,1);$

\emph{i) Algebras from $F$ are isomorphic to} $%
L(0,0,0,0,0,0).$
\end{pr}

\begin{theor}
\emph{Let $L$ be non-Lie complex filiform Leibniz algebra in
$SLeib_{8}$. \ Then it is isomorphic to one of the following
pairwise non-isomorphic Leibniz algebras:}

 $1)\ \ L(1,0,\lambda_{1},\lambda_{1},\lambda_{2},\lambda_{3}):$

$\qquad \qquad  e_{0}e_{0}=e_{2},\ \ e_{i}e_{0}=e_{i+1,}\ \ 2\leq
i\leq 6,\ \ e_{0}e_{1}=e_{3}+\lambda _{1}e_{5}+\lambda
_{1}e_{6}+\lambda  _{2}e_{7},\ \ e_{1}e_{1}=\lambda _{3}e_{7},\ $

$\qquad \qquad  e_{2}e_{1}=e_{4}+\lambda  _{1}e_{6}+\lambda
_{1}e_{7},\ \ e_{3}e_{1}=e_{5}+\lambda _{1}e_{7},\ \
e_{4}e_{1}=e_{6},$

$\qquad \qquad  e_{5}e_{1}=e_{7},\ \ \lambda  _{1},\lambda
_{2},\lambda
_{3}\in C.$\\

$2)\ \ L(1,0,\lambda _{1},0,\lambda _{2},\lambda _{1}):$

$\qquad \qquad e_{0}e_{0}=e_{2},\ \ e_{i}e_{0}=e_{i+1,}\ \ 2\leq
i\leq 6,\ \ e_{0}e_{1}=e_{3}+\lambda _{1}e_{5}+\lambda _{2}e_{7},\
\ e_{1}e_{1}=\lambda  _{1}e_{7},$

$\qquad \qquad e_{2}e_{1}=e_{4}+\lambda _{1}e_{6},\ \
e_{3}e_{1}=e_{5}+\lambda_{1}e_{7},\ \ e_{4}e_{1}=e_{6},\ \
e_{5}e_{1}=e_{7},\ \ \lambda _{1},\lambda
_{2}\in C.$\\

$3)\ \  L(1,0,1,0,\lambda  ,0):$

$\qquad \qquad  e_{0}e_{0}=e_{2},\ \ e_{i}e_{0}=e_{i+1,}\ \ 2\leq
i\leq 6,\ \ e_{0}e_{1}=e_{3}+e_{5}+\lambda e_{7},\ \
e_{2}e_{1}=e_{4}+e_{6},\ \ e_{3}e_{1}=e_{5}+e_{7},$

$\qquad \qquad  e_{4}e_{1}=e_{6},\ \ e_{5}e_{1}=e_{7},\ \ \lambda
\in
C.$\\

$4)\ \  L(1,0,0,\lambda_{1},\lambda_{1},\lambda_{2}):$

 $ \qquad \qquad e_{0}e_{0}=e_{2},\ \ e_{i}e_{0}=e_{i+1,}\ \
2\leq i\leq 6,\ \ e_{0}e_{1}=e_{3}+\lambda
_{1}e_{6}+\lambda_{1}e_{7},\ \ e_{1}e_{1}=\lambda_{2}e_{7},$

$\qquad \qquad e_{2}e_{1}=e_{4}+\lambda_{1}e_{7},\ \
e_{3}e_{1}=e_{5},\ \ e_{4}e_{1}=e_{6},\ \
e_{5}e_{1}=e_{7},\ \ \lambda_{1},\lambda_{2}\in C.$\\

$5)\ \  L(1,0,0,1,0,\lambda ):$

$\qquad  \qquad e_{0}e_{0}=e_{2},\ \ e_{i}e_{0}=e_{i+1,}\ \ 2\leq
i\leq 6,\ \ e_{0}e_{1}=e_{3}+e_{6},\ \ e_{1}e_{1}=\lambda e_{7},\
\ e_{2}e_{1}=e_{4}+e_{7},$

$\qquad  \qquad  e_{3}e_{1}=e_{5},\ \ e_{4}e_{1}=e_{6},\ \
e_{5}e_{1}=e_{7},\ \ \lambda \in
C.$\\

$6)\ \ L(1,0,0,0,\lambda,\lambda):$

$\qquad  \qquad e_{0}e_{0}=e_{2},\ \ e_{i}e_{0}=e_{i+1,}\ \ 2\leq
i\leq 6,\ \ e_{0}e_{1}=e_{3}+\lambda e_{7},\ \ e_{1}e_{1}=\lambda
e_{7},\ \  e_{2}e_{1}=e_{4},\ \ e_{3}e_{1}=e_{5},$

$ \qquad \qquad e_{4}e_{1}=e_{6},\ \ e_{5}e_{1}=e_{7},\ \ \lambda
\in
C. $\\

 $7)\ \  L(0,1,1,0,\lambda _{1},\lambda_{2}):$

$\qquad  \qquad e_{0}e_{0}=e_{2},\ \ e_{i}e_{0}=e_{i+1,}\ \ 2\leq
i\leq 6,\ \ e_{0}e_{1}=e_{4}+e_{5}+\lambda _{1}e_{7},\ \
e_{1}e_{1}=\lambda _{2}e_{7},\ \ e_{2}e_{1}=e_{5}+e_{6},$

$\qquad \qquad e_{3}e_{1}=e_{6}+e_{7},\ \ e_{4}e_{1}=e_{7},\ \
\lambda _{1},\lambda _{2}\in C.$\\

 $8)\ \ L(0,1,1,-1,\lambda _{1},\lambda _{2}):$

$\qquad \qquad e_{0}e_{0}=e_{2},\ \ e_{i}e_{0}=e_{i+1,}\ \
2\leq i\leq 6,\ \ e_{0}e_{1}=e_{4}+e_{5}-e_{6}+%
\lambda_{1}e_{7},\ \ e_{1}e_{1}=\lambda_{2}e_{7}, $

$\qquad \qquad e_{2}e_{1}=e_{5}+e_{6}-e_{7},\ \
e_{3}e_{1}=e_{6}+e_{7}, \ \ e_{4}e_{1}=e_{7},\ \ \lambda _{1},%
\lambda_{2}\in C.$\\

 $9)\ \ L(0,1,0,0,\lambda ,1):$

$ \qquad \qquad e_{0}e_{0}=e_{2},\ \ e_{i}e_{0}=e_{i+1,}\ \ 2\leq
i\leq 6,\ \ e_{0}e_{1}=e_{4}+\lambda e_{7},\ \ e_{1}e_{1}=e_{7},\
\  e_{2}e_{1}=e_{5},\ \ e_{3}e_{1}=e_{6},$

 $ \qquad \qquad e_{4}e_{1}=e_{7},\ \ \lambda\in C.$\\

 $10)\ \ L(0,1,0,-1,\lambda,1):$

 $ \qquad \qquad e_{0}e_{0}=e_{2},\ \ e_{i}e_{0}=e_{i+1,}\ \
2\leq i\leq 6,\ \ e_{0}e_{1}=e_{4}-e_{6}+\lambda e_{7},\ \
e_{1}e_{1}=e_{7},\ \ e_{2}e_{1}=e_{5}-e_{7},$

$\qquad \qquad  e_{3}e_{1}=e_{6},\ \ e_{4}e_{1}=e_{7},\ \ \lambda
\in
C.$\\

$11)\ \ L(0,0,1,1,0,\lambda):$

$\qquad \qquad e_{0}e_{0}=e_{2},\ \ e_{i}e_{0}=e_{i+1,}\ \ 2\leq
i\leq 6,\ \ e_{0}e_{1}=e_{5}+e_{6},\ \ e_{1}e_{1}=\lambda e_{7},\
\ e_{2}e_{1}=e_{6}+e_{7},$

$\qquad \qquad e_{3}e_{1}=e_{7},\ \ \lambda \in C.$\\

 $12)\ \  L(0,0,1,0,\lambda,0):$

 $ \qquad \qquad e_{0}e_{0}=e_{2},\ \ e_{i}e_{0}=e_{i+1,}\ \
2\leq i\leq 6,\ \ e_{0}e_{1}=e_{5}+\lambda e_{7},$ $\
e_{2}e_{1}=e_{6},\ \ e_{3}e_{1}=e_{7},\ \ \lambda \in
C.$\\

 $13)\ \ L(0,0,0,1,0,\lambda ):$

$ \qquad \qquad e_{0}e_{0}=e_{2},\ \ e_{i}e_{0}=e_{i+1,}\ \ 2\leq
i\leq 6,\ \ e_{0}e_{1}=e_{6},\ \
e_{1}e_{1}=\lambda e_{7},\ \ e_{2}e_{1}=e_{7},\ \ \lambda\in C.$\\

$ 14)\ \  L(0,0,0,1,1,\lambda ):$

 $ \qquad \qquad e_{0}e_{0}=e_{2},\ \ e_{i}e_{0}=e_{i+1,}\ \
2\leq i\leq 6,\ \ e_{0}e_{1}=e_{6}+e_{7},\ \
e_{1}e_{1}=\lambda e_{7},\ \ e_{2}e_{1}=e_{7},\ \ \lambda\in C.$\\

$15)\ \ L(0,1,0,0,1,0):$

 $ \qquad \qquad e_{0}e_{0}=e_{2},\ \ e_{i}e_{0}=e_{i+1,}\ \
2\leq i\leq 6,\ \ e_{0}e_{1}=e_{4}+e_{7},\ \
e_{2}e_{1}=e_{5},\ \ e_{3}e_{1}=e_{6},\ \ e_{4}e_{1}=e_{7}.$\\

$16)\ \  L(0,1,0,0,0,0):$

 $ \qquad \qquad e_{0}e_{0}=e_{2},\ \ e_{i}e_{0}=e_{i+1,}\ \
2\leq i\leq 6,\ \ e_{0}e_{1}=e_{4},\ \
e_{2}e_{1}=e_{5},\ \ e_{3}e_{1}=e_{6},\ \ e_{4}e_{1}=e_{7}.$\\

 $17)\ \  L(0,1,0,-1,1,0):$

$\qquad \qquad  e_{0}e_{0}=e_{2},\ \ e_{i}e_{0}=e_{i+1,}\ \ 2\leq
i\leq 6,\ \ e_{0}e_{1}=e_{3}-e_{6}+e_{7},\ \
e_{2}e_{1}=e_{5}-e_{7},\ \ e_{3}e_{1}=e_{6},$

$\qquad \qquad  e_{4}e_{1}=e_{7}.$\\

$18)\ \ L(0,1,0,-1,0,0):$

 $ \qquad \qquad e_{0}e_{0}=e_{2},\ \ e_{i}e_{0}=e_{i+1,}\ \
2\leq i\leq 6,\ \ e_{0}e_{1}=e_{4}-e_{6},\ \
e_{2}e_{1}=e_{5}-e_{7},\ \ e_{3}e_{1}=e_{6},\ \ e_{4}e_{1}=e_{7}.$\\

$19)\ \ L(0,0,1,0,0,1):$

 $ \qquad \qquad e_{0}e_{0}=e_{2},\ \ e_{i}e_{0}=e_{i+1,}\ \
2\leq i\leq 6,\ \ e_{0}e_{1}=e_{5},\ \
e_{1}e_{1}=e_{7},\ \  e_{2}e_{1}=e_{6},\ \ e_{3}e_{1}=e_{7}.$\\

$ 20)\ \  L(0,0,0,0,0,1):$

 $ \qquad \qquad e_{0}e_{0}=e_{2},\ \ e_{i}e_{0}=e_{i+1,}\ \
2\leq i\leq 6,\ \ e_{1}e_{1}=e_{7}.$\\

$ 21)\ \  L(0,0,0,0,1,0):$

 $ \qquad \qquad e_{0}e_{0}=e_{2},\ \ e_{i}e_{0}=e_{i+1,}\ \
2\leq i\leq 6,\ \ e_{0}e_{1}=e_{7}.$\\

$ 22)\ \  L(0,0,0,0,1,1):$

 $ \qquad \qquad e_{0}e_{0}=e_{2},\ \ e_{i}e_{0}=e_{i+1,}\ \
2\leq i\leq 6,\ \ e_{0}e_{1}=e_{7},\ \
e_{1}e_{1}=e_{7}.$\\

$ 23)\ \  L(0,0,0,0,0,0):$

$ \qquad \qquad e_{0}e_{0}=e_{2},\ \ e_{i}e_{0}=e_{i+1,}\ \ 2\leq
i\leq 6.$

\end{theor}

The adapted number of isomorphism classes $N_{8}$=23.

\textbf{Conjecture.} The adapted number of isomorphism classes
$N_n$ of $n$-dimensional non-Lie complex filiform Leibniz algebras
in $SLeib_{n}$ can be found by the formula:
$$N_n=n^2-7n+15.$$

\section{\textbf{About method of classification}}

To classify $SLeib_n$ we split it into several subsets and
classify the algebras from each of these subsets separately. The
formula $n^2-7n+15$ for the adapted number of isomorphism classes
has a hypothetic character, but it is confirmed by our
computations in dimension 9 as well.

It is a slightly tedious to check for the expressions given in the
parametric family cases to be invariant by hand, but a computer
can do it very efficiently. This procedure has been implemented in
the computer programm Maple 10.

\end{document}